\documentclass{ifacconf}
\makeatletter
\let\algorithm\@undefined
\let\endalgorithm\@undefined
\makeatother
\usepackage[ruled,vlined,linesnumbered]{algorithm2e}
\usepackage{graphicx,siunitx}      % include this line if your document contains figures
\usepackage{natbib}        % required for bibliography
\usepackage{amssymb,mathtools,theorem,psfrag}
\usepackage{amsmath,color}

\newtheorem{example}{Example}
\newtheorem{proposition}{Proposition}

\newtheorem{definition}{Definition}
\newtheorem{lemma}{Lemma}
\newtheorem{problem}{Problem}
\newtheorem{remark}{Remark}

\newcommand*{\QEDA}{\hfill\ensuremath{\blacksquare}}%
\newcommand*{\QEDB}{\hfill\ensuremath{\square}}%
\DeclareMathOperator{\q}{q}
\DeclareMathOperator{\He}{He}
\DeclareMathOperator{\dom}{dom}
\DeclareMathOperator{\Sign}{sign}

\DeclareMathOperator{\trace}{trace}
\DeclareMathOperator{\maximize}{maximize}
\newcommand{\R}{\mathbb{R}}
\begin{document}
\begin{frontmatter}

\title{Anti-windup-like Compensator Synthesis for Discrete-Time Quantized Control Systems\\ (Extended Version)}
\thanks{This work has been supported in part by the ANR Labex CIMI (grant ANR-11-LABX- 0040) within the French State Program ``Investissement d'Avenir''.}
%\thanks{Emails: samer.alsamadi@mines-albi.fr, francesco.ferrante@unipg.it, tarbour@laas.fr.} 
% Title, preferably not more than 10 words. \thanks[footnoteinfo]{This work was in part supported by the ANR Labex CIMI (grant ANR-11-LABX-0040) within the French State Programme "Investissement d'Avenir".}}

\author[IMT]{Samer Alsamadi} 
\author[UNIPG]{ Francesco Ferrante} 
\author[LAAS]{Sophie Tarbouriech}

\address[IMT]{Centre de Génie Industriel, IMT Mines Albi, Albi, France. samer.alsamadi@mines-albi.fr}
\address[UNIPG]{Department of Engineering, University of Perugia, Perugia, Italy. francesco.ferrante@unipg.it}
\address[LAAS]{LAAS-CNRS, Universit\'e de Toulouse, CNRS, Toulouse, France. tarbour@laas.fr}

\begin{abstract}                % Abstract of not more than 250 words.
This paper addresses the problem of designing an anti-windup like compensator for discrete-time linear control systems with quantized input. The proposed compensator provides a correction signal proportional to the quantization error that is fed to the controller. The compensator is designed to ensure that solutions to the closed-loop systems converge in finite time into a compact set containing the origin that can be tuned by the designer.  A numerically tractable algorithm with feasibility guarantees is provided for the design of the compensator. The proposed results are illustrated on an academic example and an open-loop unstable aircraft system.
\end{abstract}

\begin{keyword} Quantized control, anti-windup like compensator design, LMIs, stability.
\end{keyword}

\end{frontmatter}
%===============================================================================
\section{Introduction}
\subsection{Background and motivation}
Recent technological enhancements have enabled the design of a new generation of control systems combining physical interactions with computational and communication abilities. The rapid spread of these technologies is related to their enormous advantages in terms of scalability, ease of maintenance, and high computational resources. This has largely impacted numerous applications such as transportation systems,  autonomous robotics, and energy delivery systems, just to mention a few. This new trend has had also a strong impact in modern control systems, which are by now  characterized by the interplay of digital controllers and/or digital instrumentation and physical systems (\cite{mur:ast:boy:bro:ste/repor2002}).

One particular problem in this context resides in the presence of quantizers in control loops, especially when the quantization affects the input of the system. Quantization is generally used to reduce data traffic load in communication channels to comply with limited bandwidth constraints. As a side effect, quantization introduces nonlinearity into the closed loop. The unwanted effects of this nonlinearity include chaotic behaviors, additional equilibrium points, and limit cycles, as stated in \cite{delchamps1990stabilizing}, \cite{cer:dep:fra/nolcos2010},  \cite{lib/book2003}, and \cite{tarbouriech2011control}. Over the last few years, researchers have mostly focused on the analysis of quantized control systems  involving continuous-time plants/controllers; see, e.g., \cite{RogerBrockett}, \cite{DanielLiberzon}, \cite{EmiliaFridman2},  \cite{ferrante2015stabilization}, and \cite{ferrante2020sensor}. On the other hand, due to the inherent relationship between digital control and quantization, the analysis of quantization in discrete-time control systems is a relevant problem in applications. This has pushed the community to address quantization in a discrete-time setting.  Notable results on stability analysis and control design of discrete-time quantized control systems can be found in \cite{BrunoPicasso}, \cite{MinyueFu-ieee09}, \cite{campos2018stabilisation}, and \cite{ich:saw:tar/ijrnc2018}. A novel approach for stability analysis of finite-level quantized discrete-time systems has been recently proposed in \cite{valmorbida2020quantization}. 
\subsection{Contributions}
The focus of this paper is on discrete-time control systems with input quantization. In this setting, we propose an approach for compensating a pre-designed closed-loop control system to reduce the effect of input quantization. The approach we pursue is inspired by the use of anti-windup compensators in saturated control systems; see \cite{tarbouriech2011stability,zaccarian2011modern}.
In particular, we augment a standard output feedback control system with an additional static compensator loop. This loops feeds back a signal that is proportional to the quantization error into the controller dynamics. We show how this simple augmentation, if suitably tuned, can dramatically reduce the effect of the input quantization on the closed-loop system response.  A similar idea has been proposed in \cite{tarbouriech2018nonstandard} in the context of control systems subject to input backlash. In the context of continuous-time quantized control systems, the design of dynamic anti-windup loops is proposed in \cite{sofrony2015anti} to ensure specific $L_2$-gain performance. 

For this class of augmented control systems, the main contribution of our paper consists of establishing sufficient conditions to ensure uniform global finite-time convergence of the closed-loop system state into a compact set containing the origin with tunable shape. Later these conditions 
are used to devise an iterative design algorithm for the compensator that is based on semidefinite programming (\emph{SDP}). Under some mild assumptions, the algorithm is guaranteed to yield a feasible solution.

The paper is organized as follows. Section~\ref{sec:pb} states the class of systems under consideration and the problem we solve. Sufficient conditions for designing the compensator are presented in Section~\ref{sec:main}. 
Section~\ref{sec:CompDes}, building upon the results in Section~\ref{sec:main}, illustrates the iterative design algorithm we present in the paper. Section~\ref{sec:ex} showcases the results on two examples borrowed from the literature. Section~\ref{sec:conclu} ends the paper with some concluding remarks. Auxiliary results and definitions are included in Appendix~\ref{app:A}.

%\section{Preliminaries}
\label{sec:prelim}
\subsection{Notation}
The symbol $\mathbb{N}$ denotes the set of nonnegative integers, $\R$ ($\R_{\geq 0}$) represents the set of real (nonnegative) numbers, $\R^n$ is the $n$-dimensional Euclidean space, $\R^n_{>0}$ is the positive (open) orthant in $\R^n$, and $\R^{n\times m}$ represents the set of $n\times m$ real matrices, and $\mathbb{B}$ is the closed unitary ball in the Euclidean norm. The symbol $\mathbb{S}_{+}^n$ ($\mathbb{S}_{++}^n$) stands for the set of $n\times n$ symmetric positive semidefinite (definite) matrices, 
 $\mathbb{D}_+^{n}$  ($\mathbb{D}_{++}^{n}$) denotes the set of $n\times n$ diagonal positive semidefinite (definite)  matrices, and $\mathbb{P}^{n}$ is the set of $n\times n$ symmetric matrices with nonnegative entries. For a vector $x\in\R^n$ (a matrix $A\in\R^{n\times m}$) $x^\top$ ($A^\top$) denotes the transpose of $x$ ($A$). 
The spectral radius of the matrix  $A$ is denoted by $\rho(A)$.
 Given $A\in\R^{n\times n}$, $\He (A)=A+A^\top$. The symbol $A\preceq 0$ ($A\prec 0$) stands for seminegative (negative) definiteness of the symmetric matrix $A$. The symbol $\star$ stands for symmetric blocks is symmetric matrices. Given a symmetric matrix $A$, 
$\lambda_{\max}(A)$ and $\lambda_{\min}(A)$ stand, respectively, for the largest and smallest eigenvalue of $A$. Given $A\in\R^{n\times n}$, the notation $\mathcal{E}(A)=\{x\in\R^n\colon x^\top A x\leq 1 \}$ is used. The function $\Sign\colon\R\to\{-1, 1\}$ is defined for all $x\in\R$ as follows: $\Sign(x)=1$ if $x\geq 0$ and $-1$ otherwise. The symbol $\vert x\vert_{S}\coloneqq\displaystyle\inf_{y\in S}\vert x-y\vert$ denotes the distance of the point $x\in\R^n$ to the nonempty set $S\subset \R^n$.
The symbol $\lfloor x\rfloor$ indicates the floor of the real number $x$. The symbol $L_V(c)$ denotes the $c$-sublevel set of the function $V$, i.e., $L_V(c)\coloneqq\{x\in\dom V\colon V(x)\leq c\}$. 
\section{Problem Statement}
\label{sec:pb}
We consider the following discrete-time plant:
\begin{equation}\label{1}
    \begin{cases}
      x_p^{+}=A_px_p+B_p\q_\Theta(u_p)\\
      y_p=C_px_p
    \end{cases}  
\end{equation}
where $x_p\in\mathbb{R}^{n_p}$, $u_p\in\mathbb{R}^{n_u}$, $y_p\in\mathbb{R}^{n_y}$ are, respectively, the plant state, control input, and measured output. Matrices $A_p$, $B_p$, and $C_p$ are real and assumed to be known. The function $\q_\Theta$ is the so-called \emph{uniform quantizer}, which is defined next:
$$
u\mapsto q_\Theta(u_p)=(q_{\theta_1}(u_{p1}), q_{\theta_2}(u_{p2}), \dots, q_{\theta_{n_u}}(u_{p n_u}))
$$
where $\Theta=(\theta_1, \theta_2,\dots, \theta_{n_u})\in\R_{>0}^{n_u}$ represents 
the vector of quantization levels of each channel and, for all $u\in\R$, $\theta\in\R_{>0}$ 
$$
q_{\theta}(u)\coloneqq \theta  \Sign(u)\Bigl\lfloor{\frac{\vert u\vert}{\theta}\Bigr\rfloor}.
$$

The plant is controlled via the following dynamic output feedback controller:
\begin{equation}\label{2}
   \left\{ \begin{aligned}
      &x_c^{+}=A_cx_c+B_c u_c+v\\
      &y_c=C_cx_c+D_c u_c
    \end{aligned}  \right.     
\end{equation}
where $x_c\in\R^{n_c}$, $y_c\in\R^{n_u}$, and $u_c\in\R^{n_y}$ are, respectively, the controller state, input and output. The parameters $A_c$, $B_c$, $C_c$, and $D_c$ are real matrices of adequate dimensions defining the controller dynamics. The signal $v \in\R^{n_c}$ is an additional input to be designed to mitigate the effect of input quantization. This signal is reminiscent of an anti-windup correction in saturated feedback control systems \cite{zaccarian2011modern}. The use of anti-windup-like schemes in quantized control systems and systems subject to backlash have been investigated, respectively, in \cite{sofrony2015anti} and \cite{tarbouriech2018nonstandard}. Inspired by the constructions in \cite{tarbouriech2018nonstandard}, the signals $v$ is selected as follows:
\begin{equation}
\begin{aligned}
\label{eq:AWlaws}
&v=E(\q_\Theta(y_c)-y_c) = E \psi_\Theta(y_c)
\end{aligned}
\end{equation}
where $E\in\R^{n_c\times n_u}$ is a gain to be designed. Our goal is to design $E$ to reduce the effect of input quantization on the closed-loop system.  

The interconnection of the plant \eqref{1} and the controller \eqref{2}-\eqref{eq:AWlaws} is obtained by setting $u_c=y_p$ and 
$u_p=y_c$. Thus, by taking as a state $x=(x_p, x_c)$, the closed-loop system reads:

\begin{equation}\label{4}
    x^{+}=
    \underbrace{\begin{bmatrix}
       A_p+B_pD_cC_p & B_pC_c\\[0.3em]
       B_cC_p & A_c\\[0.3em]
     \end{bmatrix}}_{A_{CL}}
     x+
     \underbrace{\begin{bmatrix}
       B_p\\[0.3em]
       E\\[0.3em]
     \end{bmatrix}}_{B_{CL} + R E}
     \psi_\Theta(Hx)
\end{equation}
where $H\coloneqq \begin{bmatrix} D_cC_p & C_c\end{bmatrix}$ and $R\coloneqq\begin{bmatrix}0_{n_p\times n_c}\\I_{n_c}\end{bmatrix}$.
The general problem we address in this paper can be formalized as follows:
\begin{problem}
\label{pb:Design}
Given plant \eqref{1} and controller \eqref{2} parameters, and a closed set $\mathcal{U}$ containing the origin,
design a gain $E$ such that there exists a compact set $\mathcal{S}\subset \mathcal{U}$, containing the origin, that is uniformly globally finite-time attractive 
\emph{UGFTA} (see Definition~\ref{def:UGFTA} in Appendix~\ref{app:A})  for the closed-loop system \eqref{4}. 
\hfill $\diamond$
\end{problem}

\begin{remark}
The formulation of Problem~\ref{pb:Design} ensures that the state of the closed-loop system is bounded due to $\mathcal{S}$ being compact. The set $\mathcal{U}$, which is not necessarily compact, is introduced to enable the designer to shape the response of the system for large times. For example, $\mathcal{U}$ can be selected to ensure that the plant state converges close to zero. 
\hfill $\diamond$
\end{remark}

\section{Main results}
\label{sec:main}
To address Problem~\ref{pb:Design}, we make use of the following sector conditions originally introduced in \cite{ferrante2015stabilization} and later extended in \cite{ferrante2020sensor} for the mapping $\psi_\Theta$. 
\begin{lemma}
\label{lemm:Sector}
Let $S_1, S_2\in\mathbb{D}_{++}^{n_u}$. Then, for all $u\in\mathbb{R}^{n_u}$ the following inequalities hold:
$$
\begin{aligned}
&\psi_\Theta^\top (u)S_1\psi_\Theta(u)-\Theta ^\top S_1\Theta\leq0\\
&\psi_{\Theta}^\top(u)S_2(\psi_\Theta(u)+u)\leq 0.
\end{aligned}
$$
\hfill$\diamond$
\end{lemma}

We assume that the closed set $\mathcal{U}$ introduced in Problem~\ref{pb:Design} is defined as follows:
$$
\mathcal{U}=\{x\in\R^{n_p+n_c}\colon x^\top U x\leq 1\}
$$
where $U\in\mathbb{S}_{+}^{n_p+n_c}$ is given.

The following result provides sufficient conditions for the solution to Problem~\ref{pb:Design}.
\begin{proposition}
If there exist $P\in\mathbb{S}_{++}^{n_p+n_c}$, $S_1$, $S_2$ $\in\mathbb{D}_{+}^{n_u}$, $E\in\mathbb{R}^{n_c\times n_u}$, and $\tau\in (0, 1)$ such that the following conditions hold:
\begin{subequations}
\begin{align}
&\begin{bmatrix}
       (\tau-1)P & -H^\top S_2 & A^\top_{CL}P\\[0.3em]
       \star & -S_1-2S_2%-\textcolor{red}{\He(S_2E)} 
       & (B_{CL}+RE)^\top P\\[0.3em]
       \star & \star & -P\\[0.3em]
     \end{bmatrix}\prec 0\label{eq:MIuq}\\
          \label{eq:MI_uq2}
     &\Theta^\top S_1\Theta-\tau\leq 0\\
     \label{eq:inclusion}
     &U-P\preceq 0.
\end{align}
\label{eq:MainConditions}
\end{subequations}
Then, $\mathcal{S}=\mathcal{E}(P)$
is included in $\mathcal{U}$ and is UGFTA for the closed-loop system \eqref{4}.
\end{proposition}
\begin{pf}
The inclusion $\mathcal{S}=\mathcal{E}(P)\subset\mathcal{U}$ follows directly from \eqref{eq:inclusion}. The remainder of the proof is based on Proposition~\ref{prop:LyapUGFTA}. In particular, we show that the function $V(x)=x^\top Px$ satisfies all the assumptions in Proposition~\ref{prop:LyapUGFTA} with $c=1$. From Schur complement and a simple congruence transformation, the satisfaction of \eqref{eq:MIuq} implies that 
\begin{equation}
\label{eq:Mmatrix}
\begin{aligned}
M\coloneqq &\begin{bmatrix}
       (\tau-1)P & -H^\top S_2 \\[0.3em]
       \star & -S_1-2S_2%-\He(S_2E)
       \end{bmatrix} \\
       &\quad\quad\quad\quad+\begin{bmatrix}
       A_{CL}^\top\\
       (B_{CL}+RE)^\top
       \end{bmatrix}P\begin{bmatrix}
       A_{CL}^\top\\
       (B_{CL}+RE)^\top
       \end{bmatrix}^\top\!\!\!\prec 0
       \end{aligned}
\end{equation}
In particular, observe that, from simple calculations, one has
$$
\begin{aligned}
\begin{bmatrix}
   x\\
   \psi_\Theta
\end{bmatrix}^\top M\begin{bmatrix}
   x\\
  \psi_\Theta
\end{bmatrix}=&\Delta V(x)+\tau x^\top Px-\psi_\Theta^\top S_1\psi_\Theta\\
&-2\psi_\Theta^\top S_2\psi_\Theta-2\psi_\Theta^\top S_2Hx%\\
%&\textcolor{red}{-2\psi_\Theta^\top S_2E\psi_\Theta}
\end{aligned}
$$
where, for all $x\in\R^{n_p+n_c}$
$$
\begin{aligned}
\Delta V(x)\coloneqq &V(A_{CL}x+(B_{CL}+ R E)\psi_\Theta)-V(x)
\end{aligned}
$$
and the shorthand notation $\psi_\Theta=\psi_\Theta(Hx)$ is used.
Using \eqref{eq:MI_uq2} and Lemma~\ref{lemm:Sector}, one has, for all $x\in\R^{n_p+n_c}$:
 $$
\Delta V(x)+\tau x^\top P x\leq \begin{bmatrix}
   x\\
   \psi_\Theta
\end{bmatrix}^\top M\begin{bmatrix}
   x\\
   \psi_\Theta
\end{bmatrix}-\Theta^\top S_1\Theta+\tau.
 $$
 The latter, for all $x\in\R^{n_p+n_c}$, yields
$$
 \Delta V(x)+\tau (V(x)-1)\leq   \begin{bmatrix}
   x\\
   \psi_\Theta
\end{bmatrix}^\top M\begin{bmatrix}
   x\\
   \psi_\Theta
\end{bmatrix}
$$
which, by using \eqref{eq:Mmatrix}, gives
\begin{equation}
\label{eq:DeltaVSproc}
 \Delta V(x)+\tau (V(x)-1)\leq  -\varrho V(x),
\end{equation}
for some small $\varrho\in (0, 1)$. Hence, by S-procedure, the last inequality implies that, for all $x\in\overline{\R^{n_p+n_c}\setminus L_V(1)}$,
$$
 V(A_{CL}x+(B_{CL}+ R E)\psi_\Theta)\leq e^{-\mu}V(x),
$$
with $\mu\coloneqq \ln(1-\varrho)$, which corresponds to \eqref{eq:DeltaVMain}. To conclude the proof, we show that \eqref{eq:MIuq} and \eqref{eq:MI_uq2} imply \eqref{eq:InvarianceMain}. From \eqref{eq:DeltaVSproc}, for all $x\in\R^{n_p+n_c}$
$$
V(A_{CL}x+(B_{CL}+ RE)\psi_\Theta)-1+(1-\tau)(1-V(x))\leq  0
$$
which, due to $\tau\in (0, 1)$, by using S-procedure, yields
$$
V(A_{CL}x+(B_{CL}+ RE)\psi_\Theta)-1\leq 0\quad\forall x\in L_V(1),
$$
thereby giving \eqref{eq:InvarianceMain}. This establishes the result.\QEDA
\end{pf}
The result given next, which plays a relevant role in the construction of the design algorithm presented in Section~\ref{sec:CompDes}, shows that \eqref{eq:MIuq} and \eqref{eq:MI_uq2} are always feasible 
 as long as the quantization-free uncompensated ($E=0$) system closed-loop system \eqref{4} is asymptotically stable.  
 \medskip
 
\begin{proposition}
\label{prop:feasible}
If $\rho(A_{CL})<1$, then \eqref{eq:MIuq} and \eqref{eq:MI_uq2} are feasible with $E=0$ and any $\tau\in(0, 1)$  such that $\sqrt{1-\tau}>\rho(A_{CL})$. \QEDB
\end{proposition}
\begin{pf}
To show the result, we prove that \eqref{eq:MIuq} and \eqref{eq:MI_uq2} are feasible with $E=0$ for some $P\in\mathbb{S}^{n_p+n_c}_{++}$, $S_1\in\mathbb{D}^{n_u}_{++}$, and $S_2=0$. To this end,
let $\tau\in (0, 1)$ such that 
     $$
     \sqrt{1-\tau}>\rho(A_{CL}),
     $$
     this is always possible due to $\rho(A_{CL})<1$. For this selection of $\tau$, select $S_1\in\mathbb{D}^{n_u}_{++}$ such that \eqref{eq:MI_uq2} holds. Let $Q\in\mathbb{S}^{n_p+n_c}_{++}$ be any solution to the following matrix inequality:
$$
(\tau-1)Q+A^\top_{CL}QA_{CL}\prec 0
$$
which is solvable due to the selection of $\tau$ above. The latter, thanks to Schur complement lemma and a simple congruence tranformation, is equivalent to:
     \begin{equation}
 \label{eq:Elimin}
\begin{bmatrix}
       (\tau-1)Q & A^\top_{CL}Q\\[0.3em]
       \star & -Q\\[0.3em]
     \end{bmatrix}\prec 0.
\end{equation}
At this stage, observe that \eqref{eq:Elimin} can be equivalently rewritten as:
     \begin{equation}
 \label{eq:Elimin2}
\begin{bmatrix}
I&0&0\\[0.3em]
0&0&I
     \end{bmatrix}
\begin{bmatrix}
       (\tau-1)Q & 0& A^\top_{CL}Q\\[0.3em]
       \star & 0& B_{CL}^\top Q\\[0.3em]
       \star & \star & -Q\\[0.3em]
     \end{bmatrix}\begin{bmatrix}
I&0\\[0.3em]
0&0\\[0.3em]
0&I
\end{bmatrix}\prec 0.
\end{equation}
Thus, from the \emph{Projection Lemma} \cite[Lemma 3.1]{gahinet1994linear}, the satisfaction of \eqref{eq:Elimin2} implies that there exists $X\in\R^{n_u\times n_u}$ such that
%\begin{equation}
%\label{eq:Elim2}
%\begin{bmatrix}
%       (\tau-1)P & 0& A^\top_{CL}P\\[0.3em]
%       \star & 0& B_{CL}^\top P\\[0.3em]
%       \star & \star & -P\\[0.3em]
%     \end{bmatrix}+\He\left(\begin{bmatrix}
%     0\\I\\0
%     \end{bmatrix}X
%     \begin{bmatrix}
%     0&I &0
%     \end{bmatrix}
%     \right)\prec 0
% \end{equation}
\begin{equation}
\label{eq:Elim2}
\begin{bmatrix}
       (\tau-1)Q & 0& A^\top_{CL}Q\\[0.3em]
       \star & \He(X)& B_{CL}^\top Q\\[0.3em]
       \star & \star & -Q\\[0.3em]
     \end{bmatrix}\prec 0.
 \end{equation}
 Notice that, necessarily, $\He(X)\prec 0$. Let 
 $$
 \chi\leq\frac{\lambda_{\min}(S_1)}{\lambda_{\max}(-\He(X))},
 $$
observe that $\chi>0$, due to $S_1\succ 0$ and $\He(X)\prec 0$. Therefore, by setting 
$$
P\coloneqq \chi Q, \quad Z\coloneqq \chi X,
$$
from \eqref{eq:Elim2}, one has
\begin{equation}
\label{eq:Elim3}
\begin{bmatrix}
       (\tau-1)P & 0& A^\top_{CL}P\\[0.3em]
       \star & \He(Z)& B_{CL}^\top P\\[0.3em]
       \star & \star & -P\\[0.3em]
     \end{bmatrix}\prec 0.
\end{equation}
At this stage, notice that by construction, $-S_1-\He(Z)\preceq~0$, therefore \eqref{eq:Elim3} implies
$$
\begin{bmatrix}
       (\tau-1)P & 0& A^\top_{CL}P\\[0.3em]
       \star &-S_1& B_{CL}^\top P\\[0.3em]
       \star & \star & -P\\[0.3em]
     \end{bmatrix}\prec 0
$$
which corresponds to \eqref{eq:MIuq} with $S_2=0$ and $E=0$. Hence, the result is established. \QEDA
\end{pf}
\section{Compensator design and optimization aspects}
\label{sec:CompDes}
In this section, we show how the results proposed in this paper can be used to devise a computationally affordable design algorithm for the compensator gain $E$ based on SDP tools. 
\subsection{Optimization}
The formulation of Problem~\ref{pb:Design} is based on a preassigned set $\mathcal{U}$ and hence on a specific selection of the matrix $U$. On the other hand, in practice the matrix $U$ can be left as an extra degree of freedom and the design of the compensator gain $E$ recast as the following optimization 
problem:
\begin{equation}
\label{eq:opti}
\begin{aligned}
& \underset{P, E, S_1, S_2, \tau, U}{\maximize}
& & \omega(U)\\
& \text{subject to}
& & \eqref{eq:MainConditions}, h(U)=0\\
&&& U\in\mathbb{S}^{n_p+n_c}_{+}, P\in\mathbb{S}^{n_p+n_c}_{++}, S_1, S_2\in\mathbb{D}^{n_u}_{+},  
\end{aligned}
\end{equation}
where $\omega\colon\mathbb{S}^{n_p+n_c}_{+}\to\R_{\geq 0}$ associates to $U$ a suitable ``size'' and $h\colon\mathbb{S}^{n_p+n_c}_{+}\to \R^{(n_p+n_c)\times(n_p+n_c)}$ prescribes structural contraints on 
$U$.  As an example, if the primary objective is to keep the plant state $x_p$ as close as possible to zero, then a possible selection of the functions $\omega$ and $h$ is as follows:
%$$
%U=\begin{bmatrix}
%U_{11}&0\\
%0&0
%\end{bmatrix}
%$$
%where $U_{11}\in\S_{+}^{n_p+n_c}$ and the
%
%Thus, the functions $\omega$ and $h$ in this case 
$$
h(U)=\begin{bmatrix}
0_{n_p\times n_p}&U_{12}\\
\star &U_{22}
\end{bmatrix}, \omega(U)=\trace(U)
$$
where $U_{1, 2}\in\R^{n_p\times n_u}$ and $U_{2, 2}\in\R^{n_c\times n_c}$ are the corresponding blocks of the matrix $U$.  
Indeed, the selection of $h$ implies that
$$
\mathcal{U}=\{(x_p, x_c)\in\mathbb{R}^{n_p+n_c}\colon x_p^\top U_{1, 1} x_p\leq 1\}
$$
for some $U_{1,1}\in\mathbb{S}^{n_p}_{+}$, while the selection of $\omega$ ensures that the 
size of the ellipsoidal set $\mathcal{E}(U_{1, 1})$ is minimized. 
\begin{remark}
Typically, the constraint induced by the function $h$ can be eliminated by suitably structuring the matrix $U$. Therefore, henceforth such a constraint will be dropped in optimization problem \eqref{eq:opti}. 
\end{remark}
\subsection{SDP-based compensator design}
Although the function $\omega$ can be generally selected as a linear function, the fact that \eqref{eq:MIuq} is bilinear in the decision variables $E$, $\tau$, and $P$ renders 
optimization problem~\eqref{eq:opti} numerically intractable. To overcome this drawback, next we 
show how optimization problem~\eqref{eq:opti} can be (suboptimally) solved via a sequence of semidefinite programs, i.e., optimization problems with linear objective over linear matrix inequality constraints. To this end, we rely on the convex-concave decomposition approach proposed in \cite{dinh2011combining}. In a nutshell, such an approach consists of expressing bilinear terms via a convex-concave decomposition (this is always possible). As second step, provided an initial feasible point is available, the concave terms are linearized around the given feasible point. Later, the resulting linearized (SDP) problem is solved and the solution obtained is used to linearize again the original concave terms. This basically leads to a sequence of SDP problems that can be solved iteratively. This approach enjoys two interesting properties that makes it appealing to devise an iterative design algorithm. If the initial point is feasible for the original problem, then the algorithm never terminates due to infeasibility; the initial point provides always a feasible solution. 
Another key feature of this approach is that, since any psd-concave function is upper bounded by its linearization (see Lemma~\ref{lemm:DiffConv}), feasible solutions to the linearized problem are feasible for the original problem. 

To deploy this approach for the solution to optimization problem~\eqref{eq:opti}, as first step we rewrite \eqref{eq:MIuq} in the following equivalent linear-bilinear decomposed form:
$$
\underbrace{\begin{bmatrix}
      -P & -H^\top S_2 & A^\top_{CL}P\\[0.3em]
       \star & -S_1-2S_2& B_{CL}^\top P\\[0.3em]
       \star & \star & -P\\[0.3em]
     \end{bmatrix}}_{\mathcal{L}(P, S_1, S_2)}\!\!+\He\!\!\left(\!\underbrace{\begin{bmatrix}
\frac{\tau}{2}I&0\\[0.3em]
0&E^\top R^\top\\[0.3em]
0&0
 \end{bmatrix}}_{\mathcal{X}^\top(\tau, E)}\!\!\underbrace{\begin{bmatrix}
P&0&0\\[0.3em]
0&0&P
 \end{bmatrix}}_{\mathcal{Y}(P)}\right)\!\!\prec\! 0.
$$
The latter, dropping the dependency on the decision variables, can be equivalently rewritten in the following psd convex-concave decomposed form:
$$
\mathcal{L}+\mathcal{X}^\top\mathcal{X}+\mathcal{Y}^\top\mathcal{Y}-(\mathcal{X}-\mathcal{Y})^\top (\mathcal{X}-\mathcal{Y})\prec 0
$$
which, by Schur complement's lemma, is equivalent to:
\begin{equation}
\begin{bmatrix}
\mathcal{L}-\mathcal{X}^\top\mathcal{X}-\mathcal{Y}^\top\mathcal{Y}+\He(\mathcal{X}^\top\mathcal{Y})&\mathcal{X}^\top&\mathcal{Y}^\top\\
\star&-I&0\\
\star&\star&-I
\end{bmatrix}\prec 0.
\label{eq:ConvConc}
\end{equation}

The last step consists of linearizing constraint \eqref{eq:ConvConc}. To this end, we first compute the differential of the psd-concave term in  \eqref{eq:ConvConc}, i.e.
$$
\begin{aligned}
&\mathcal{Q}\coloneqq-\mathcal{X}^\top\mathcal{X}-\mathcal{Y}^\top\mathcal{Y}+\He(\mathcal{X}^\top\mathcal{Y})=\\
&\quad\quad\quad\quad\quad\quad\begin{bmatrix}
-\frac{\tau^2}{4}I-P^2+\tau P&0 &0\\
\star&E^\top R^\top R E&E^\top R^\top P\\
\star&\star&-P^2
\end{bmatrix}.
%-\begin{bmatrix}
%P^2&0 &0\\
%\star&0&0\\
%\star&\star&P^2
%\end{bmatrix}+\begin{bmatrix}
%\frac{\tau}{2}P&0 &0\\
%\star&0&E^\top R^\top P\\
%\star&\star&0
%\end{bmatrix}
\end{aligned}
$$
More precisely, we compute the differential of the following mapping
$$
(\tau, P, E) \in \R\times\mathbb{S}_{++}^{n_p+n_c}\times\R^{n_c\times n_u} %\ni (\tau, P, E)
\mapsto\mathcal{Q}(\tau, P, E),
$$
at $(\tau, P, E)\in\R\times\mathbb{S}_{++}^{n_p+n_c}\times\R^{n_c\times n_u}$. This gives:
$$
\begin{aligned}
h\mapsto (D\mathcal{Q}(\tau, P, E))h=&\begin{bmatrix}
-\frac{\tau }{2}I-P&0&0\\
\star&0&0\\
\star&0&0
\end{bmatrix}h_\tau+\\
&\begin{bmatrix}
\tau h_P+\He(Ph_P)&0&0\\
\star&0&E^\top R^\top h_P\\
\star&\star&-\He(Ph_P)
\end{bmatrix}+\\
&\begin{bmatrix}
0&0&0\\
\star&\He(E^\top  R R^\top h_E)&h_E^\top R^\top  P\\
\star&\star&0
\end{bmatrix},
\end{aligned}
$$
where the notation $h=(h_\tau, h_P, h_E)\in \R\times\mathbb{S}_{++}^{n_p+n_c}\times\R^{n_c\times n_u}$ is used. At this stage, given $(\tau_0, P_0, E_0)\in\R\times\mathbb{S}_{++}^{n_p+n_c}\times \R^{n_c\times n_u}$, the ``linear inner approximation'' of optimization problem \eqref{eq:opti} around $q^{(0)}=(\tau_0, P_0, E_0)$ reads:
$$
O^{q^{(0)}}\colon\left\{\begin{aligned}
& \underset{P, S_1, S_2, \tau, U, E}{\maximize}
& & \omega(U)\\
&\text{s.t.}
& &\!\!\!\mathcal{M}(P, S_1, S_2, E, \tau\vert P_0, E_0, \tau_0)\!\!\prec\! 0,\\
&&&\!\!\! \Theta^\top S_1\Theta-\tau\leq 0,\\
&&&\!\!\! U-P\preceq 0,\\
&&&\!\!\! U\in\mathbb{S}^{n_p+n_c}_{+}, P\in\mathbb{S}^{n_p+n_c}_{++}\\
&&&\!\!\! S_1, S_2\in\mathbb{D}^{n_u}_{+},
\end{aligned}\right.
$$
where:
$$
\begin{aligned}
\underbrace{\begin{bmatrix}
\mathcal{R}(P, S_1, S_2, E, \tau\vert P_0, E_0, \tau_0)&\mathcal{W}(\tau, P, E)\\
\star&-I
\end{bmatrix}}_{\mathcal{M}(P, S_1, S_2, E, \tau\vert P_0, E_0, \tau_0)}
\end{aligned}
$$
and
$$
\begin{aligned}
&\mathcal{R}(P, S_1, S_2, E, \tau\vert P_0, E_0, \tau_0)\coloneqq\mathcal{L}(P, S_1, S_2, P)\\
&+\mathcal{Q}(\tau_0, P_0, E_0)+(D\mathcal{Q}(\tau_0, P_0, E_0))(\tau-\tau_0, P-P_0, E-E_0)\\
&\mathcal{W}(\tau, P, E)\coloneqq\begin{bmatrix}\mathcal{X}^\top(\tau, E)&\mathcal{Y}^\top(P)\end{bmatrix}.
\end{aligned}
$$

As mentioned earlier, the applicability of the convex-concave decomposition approach in \cite{dinh2011combining} requires the knowledge of an initial feasible solution to optimization problem~\eqref{eq:opti}.  To this end, we make use of Proposition~\ref{prop:feasible} and select the feasible initial point as the solution to the following optimization problem:

\begin{equation}
\label{eq:opti2}
\begin{aligned}
& \underset{P, S_1, S_2, \tau, U}{\maximize}
& & \omega(U)\\
& \text{subject to}
&& \begin{bmatrix}
       (\tau-1)P & -H^\top S_2 & A^\top_{CL}P\\
       \star & -S_1-2S_2& B_{CL}^\top P\\
       \star & \star & -P\\
     \end{bmatrix}\prec 0\\
   & & &\Theta^\top S_1\Theta-\tau\leq 0\\
   & &  &U-P\preceq 0\\
   &&& U\in\mathbb{S}^{n_p+n_c}_{+}, P\in\mathbb{S}^{n_p+n_c}_{++}, S_1, S_2\in\mathbb{D}^{n_u}_{+}.
\end{aligned}
\end{equation}
Proposition~\ref{prop:feasible} ensures that \eqref{eq:opti2} is always feasible provided that $\rho(A_{CL})<1$. Moreover, since whenever $\tau$ is fixed (and $\omega$ is linear) \eqref{eq:opti2} is an SDP program, a feasible solution to \eqref{eq:opti2} can be easily computed by performing a line search on the variable $\tau$ in the interval $(0, 1)$. 

Based on the steps presented so far, our approach to solve optimization problem~\eqref{eq:opti} is summarized in Algorithm~\ref{alg:design}.
\IncMargin{2em}
\begin{algorithm}
\caption{Optimal compensator synthesis}
\label{alg:design}
\Indm
\KwIn{Matrices $A_{CL}, B_{CL}$, quantization levels vector $\Theta$, a linear function $\omega\colon\mathbb{S}^{n_p+n_c}\to\R$, $k_{\max}\in\mathbb{N}_{>0}$, and $\varepsilon>0$.}
\Indp
  \BlankLine
  
\textbf{Initial solution:} Solve \eqref{eq:opti2} via a line search on $\tau\in(0, 1)$. Let $P_0, \tau_0$, and $U_0$ the values associated to the corresponding solution. Set $k=0$ and $E_0=0$\;
\While{$k<k_{\max}$}{
Solve SDP problem $O^{(\tau_k, P_k, E_k)}$\;
$\tau_{k+1}\longleftarrow \tau$, $P_{k+1}\longleftarrow P$, $U_{k+1}\longleftarrow U$, $E_{k+1}\longleftarrow E$\;
\SetAlgoNoLine
\If{$\vert\omega(U_{k+1})-\omega(U_k)\vert\leq \varepsilon$}{break\;}
$k\longleftarrow k+1$\;
}
\KwRet{$E$}
\end{algorithm}
\begin{remark}
An alternative approach to handle the bilinear constraint \eqref{eq:MIuq} consists of the use of a coordinate descent-type algorithm, i.e., alternatively and iteratively fixing some of the decision variables and optimizing with respect to the others; see, e.g.,
\cite{peaucelle2001efficient} and the references therein. However, numerical experiments show that this latter approach leads to more conservative results. 
\end{remark}
\section{Numerical Examples}
\label{sec:ex}
In this section, we showcase the
application of the methodology proposed in the paper in two examples. 
The first example, which is more academic, pertains to an unstable nonminimum phase plant. The second example is of practical interest and concerns a linearized model of an aircraft.
Numerical solutions to LMIs are obtained through the solver \textit{MOSEK} \cite{mosek} and coded in Matlab$^{\tiny{\textregistered}}$ via \textit{YALMIP} \cite{lofberg2004yalmip}. 

\begin{example}
\label{ex:ex1}
We consider the following nonminimum phase unstable discrete-time plant borrowed from \cite{fu2005sector}:
\begin{equation}
\label{eq:NMP}
\begin{aligned}
&x^+_p=\begin{bmatrix}
0 & 1 & 0\\ 0 & 2& 0\\ -3 & 1 & 0 
\end{bmatrix}x_p+\begin{bmatrix}
0\\1\\0
\end{bmatrix}\q_\Theta(u_p)\\
&y_p=\begin{bmatrix} -3    & 1  &   0\end{bmatrix}x
\end{aligned}
\end{equation}
We assume that $\Theta=0.5$ and that the plant is controlled by the following LQG feedback stabilizing controller:
$$
\begin{aligned}
&x_c^+=\begin{bmatrix}
-4.6 & 2.53 & 0\\ -9.2 & 3.39 & 0\\ -0.0609 & 0.0203 & 0 
\end{bmatrix}x_c+\begin{bmatrix}
-1.53\\ -3.07\\ 0.98 
\end{bmatrix}u_c\\
&y_c=\begin{bmatrix}0& -1.67 & 0\end{bmatrix}x_c.
\end{aligned}
$$
In this example, we select the following structure of the matrix $U$:
$$
U=c\begin{bmatrix}I&0\\0&0\end{bmatrix}
$$
where $c\geq 0$ is decision variable and $\omega(U)=\trace(U)$. Setting $\varepsilon=10^{-4}$, Algorithm~\ref{alg:design}  terminates in 57 iterations and returns the following gain for the compensator:
$$
E=\begin{bmatrix}
0.0379\\ 1.0645\\ 0.01
\end{bmatrix}.
$$ 

In \figurename~\ref{fig:Ex1Response} we report the evolution of the plant state, from the initial condition $x_p(0)=(1, 2, -1), x_c(0)=0$, obtained with and without the use of the compensator. The picture clearly shows that the proposed compensation strategy leads to a dramatic improvement in the plant state response. 
\begin{figure}
\psfrag{k}[1][1][1]{$j$}
\psfrag{xp}[1][1][1]{$x_p$}
\includegraphics[width=\columnwidth, trim=0.1cm 0 1cm 0.5cm, clip]{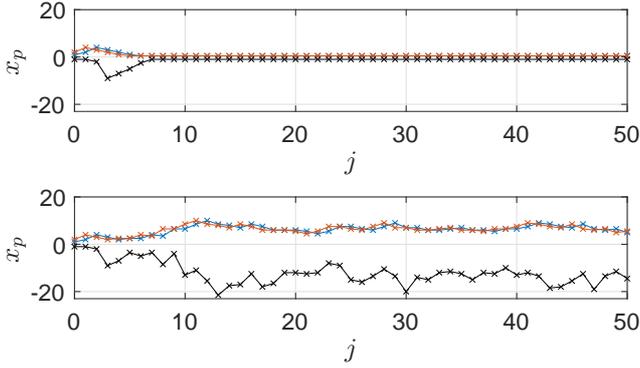}
\caption{Example~\ref{ex:ex1}. Evolution of the plant state from $x_p(0)=(1, 2, -1)$, $x_c(0)=0$. With compensation (top) and without compensation (bottom).}
\label{fig:Ex1Response}
\end{figure}
\end{example}
\begin{example}
\label{ex:ex2}
We consider the discretized\footnote{System \eqref{eq:TAFA} is obtained by performing a ZOH-discretization of the model in \cite{cristofaro2019switched} with a sampling period $T_s=0.1$.} linearized model of the short period longitudinal dynamics of TAFA (Tailless Advanced Fighter Aircraft) in \cite{cristofaro2019switched} 
\begin{equation}
\label{eq:TAFA}
\begin{aligned}
&\begin{bmatrix}
	\alpha^+\\ \psi^+\\
\end{bmatrix}=\begin{bmatrix}
0.94 & 0.087\\ 0.516 & 0.836
\end{bmatrix}\begin{bmatrix}
	\alpha\\\psi
\end{bmatrix}+\begin{bmatrix}
0.0364\\ 0.729
\end{bmatrix}u\\
&y=\begin{bmatrix}0&1\end{bmatrix}\begin{bmatrix}
	\alpha\\\psi
\end{bmatrix}.
\end{aligned}
\end{equation}
The variables $\psi$ (\SI{}{\radian\per\second}) and $\alpha$ (\SI{}{\radian}) represent, respectively, the body axis pitch rate
and the deviation of the angle of attack. The control input $u$ (\SI{}{\radian}) corresponds to the
deviation of the elevator deflection. We focus on a scenario in which the control input is quantized via a uniform quantizer with quantization level $\Theta=0.0035$, which corresponds to a quantization of  $0.2\,\, [\SI{}{\deg}]$ of the elevator deflection. The system is controlled via the following output feedback stabilizing controller:
$$
\begin{aligned}
&x_c^+=\begin{bmatrix}0.706 & -1.58\\ -4.17 & -1.88\end{bmatrix}x_c+\begin{bmatrix}
1.62\\ 1.68
\end{bmatrix}u_c\\
&y_c=\begin{bmatrix}-6.43 & -1.43\end{bmatrix}x_c.
\end{aligned}
$$
In this example, we select  $U=P$
and $\omega(U)=\trace(U)$. Setting $\varepsilon=10^{-3}$, Algorithm~\ref{alg:design} terminates in 736 iterations and returns the following gain for the compensator:
$$
E=\begin{bmatrix}
-0.0775\\
0.7222
\end{bmatrix}.
$$ 

In \figurename~\ref{fig:Ex2Response} we report a simulation of the evolution of the closed-loop plant state and of the control input from the initial condition $x_p(0)=(\frac{\pi}{180} , 0), x_c(0)=0$. To further emphasize the benefits of the proposed compensation strategy, in this simulation the compensation is activated at $j=50$ and deactivated again at $j=100$. It is interesting to notice that the use of the compensator not only leads to an improved plant state response but also to a reduced control effort. 
\begin{figure}
\psfrag{k}[1][1][1]{$j$}
\psfrag{u}[1][1][1]{$u\,\, [\SI{}{\deg}]$}
\psfrag{a}[1][1][1]{$\alpha\,\, [\SI{}{\deg}]$}
\psfrag{psi}[1][1][1]{$\psi\,\, [\SI{}{\deg\per\second}]$}
\includegraphics[width=\columnwidth, trim=0.1cm 0 1cm 0.5cm, clip]{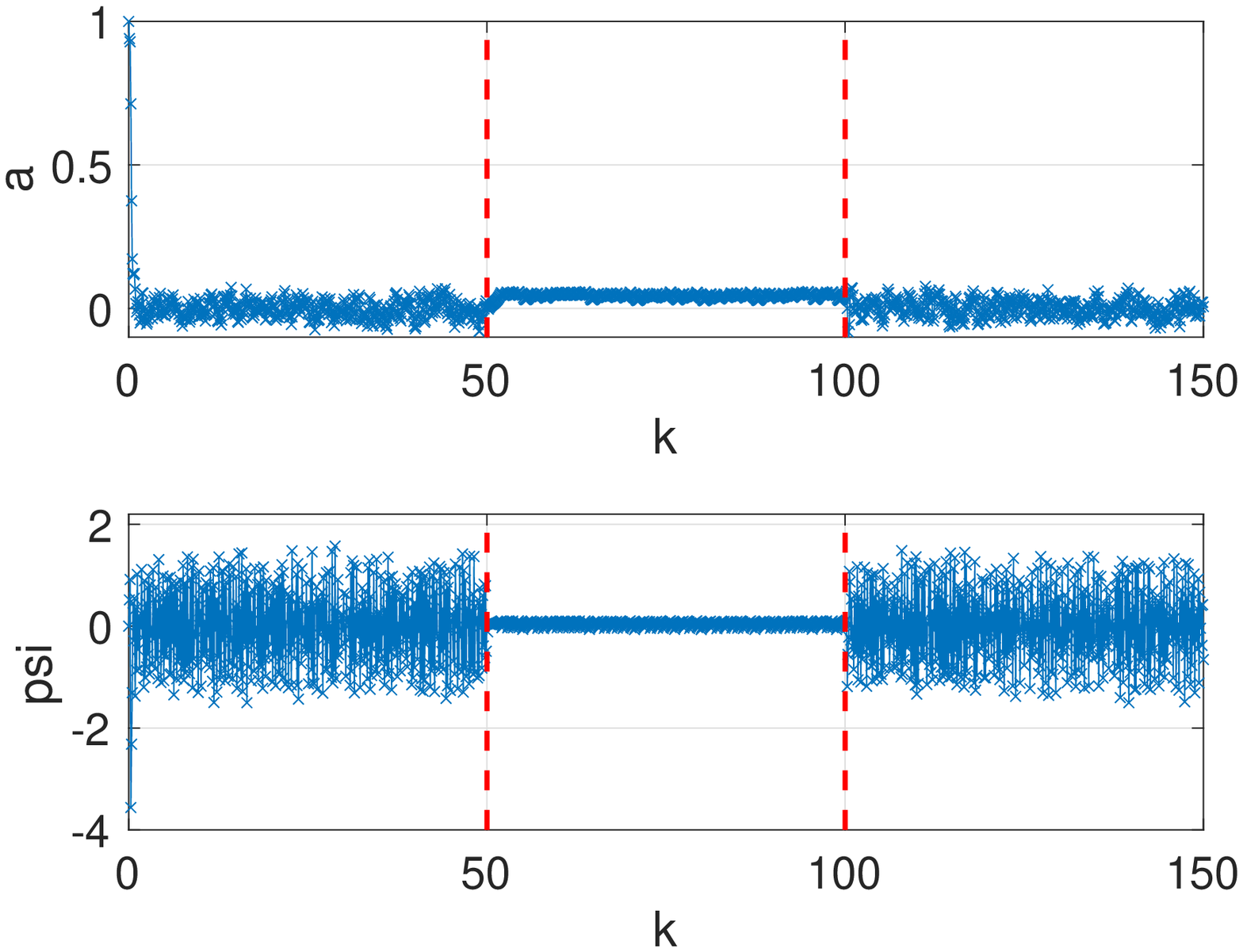}
\includegraphics[width=\columnwidth, trim=0.1cm 0 1cm 0.1cm, clip]{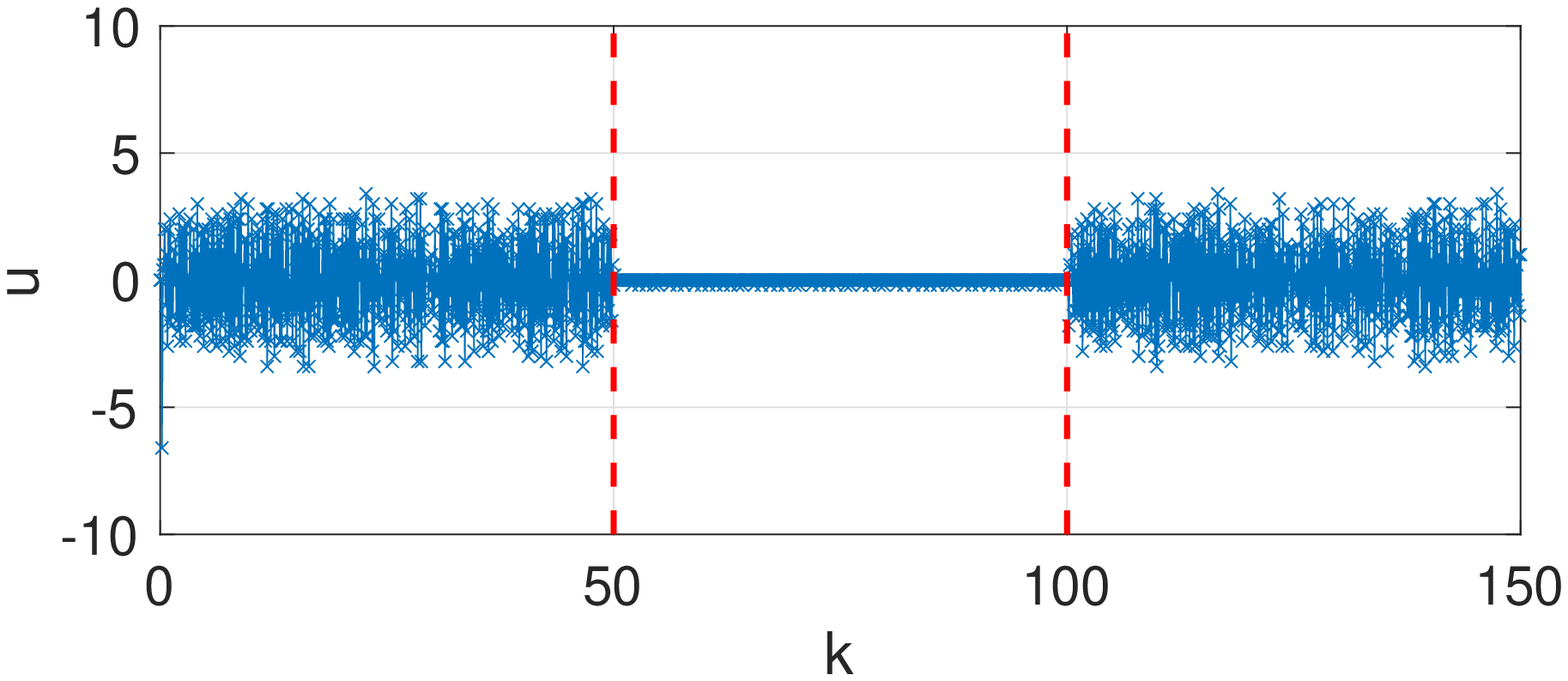}
\caption{Example~\ref{ex:ex2}. Evolution of the plant states from $x_p(0)=(\frac{\pi}{180} , 0), x_c(0)=0$. The compensation is activated at $j=50$ and deactivated at $j=100$.}
\label{fig:Ex2Response}
\end{figure}

\end{example}
\section{Conclusion}
\label{sec:conclu}
This paper addressed the design of a static anti-windup-like loop for linear closed-loop control systems subject to quantized input. The compensator was designed so as to reduce the size of the closed-loop attractor. The approach relied on matrix inequality stability conditions and result into an iterative algorithm based on a sequence of semidefinite programs for the design of the compensator gain. The effectiveness of the compensation mechanism was validated via two numerical examples.

This work opens the door to several research lines. In particular, the addition of a supplementary correction signal affecting the output of the controller is under study. This modification is particularly challenging since it leads to discontinuous algebraic loops, which cannot be addressed via existing results for saturated systems.  In addition, the extension of the proposed approach to sampled-data systems is currently part of our research. 
\bibliography{biblio}     
\appendix
\section{Auxiliary results and definitions}
\label{app:A}
\subsection{Discrete-time systems}
In this paper, we consider autonomous nonlinear discrete-time systems of the form:
\begin{equation}
\label{eq:NLsys}
x^+=g(x)
\end{equation}
where $x\in\R^n$ is the system state, $g\colon\R^n\to\R^n$, and $x^+$ stands for the value of $x$ after a jump. A function  $\phi\colon\dom\phi\rightarrow\R^n$ is a solution to \eqref{eq:NLsys} if $\dom\phi=\mathbb{N}\cap \{0, 1,\dots, J\}$ for some $J\in\mathbb{N}\cup\{\infty\}$ and for all $j\in\dom\phi$ such that $j+1\in\dom\phi$, $\phi(j+1)=g(\phi(j))$. A solution $\phi$ to \eqref{eq:NLsys} is said to be maximal if its domain cannot be extended and complete if $\dom\phi$ is unbounded. In particular, notice that maximal solutions to \eqref{eq:NLsys} are complete.    
\begin{definition}
\label{def:UGFTA}
Let $S\subset\R^n$ be closed. We say that $S$ is uniformly globally finite-time attractive (UGFTA) for \eqref{eq:NLsys} if there exists a locally bounded function $T\colon\mathbb{R}_{\geq 0}\rightarrow \mathbb{N}_{\geq 0}$ such that for any maximal solution $\phi$ to \eqref{eq:NLsys}, $j\geq T(\vert\phi(0)\vert_{S})$ implies that $\phi(j)\in S$.\hfill$\diamond$
\end{definition}
Next we provide a sufficient condition for a sublevel set of a function $V$ to be UGFTA for \eqref{eq:NLsys}.   
\begin{proposition}
\label{prop:LyapUGFTA}
Let $V\colon\R^n\rightarrow\R_{\geq 0}$ be locally bounded, $\mu<0$, $c>0$ such that $L_V(c)$ is compact, and
\begin{align}
\label{eq:DeltaVMain}
&V(g(x))\leq e^{\mu} V(x)&\forall x\in\overline{\R^n\setminus L_V(c)}\\
&V(g(x))\leq c&\forall x\in L_V(c)
\label{eq:InvarianceMain}
\end{align}
Then, $L_V(c)$ is UGFTA for \eqref{eq:NLsys}.
\end{proposition}
\begin{pf}
As a first step, notice that \eqref{eq:InvarianceMain} implies that $L_V(c)$ is forward invariant for \eqref{eq:NLsys}. Now, we show that \eqref{eq:DeltaVMain} implies that any maximal solution to \eqref{eq:NLsys} converges to $L_V(c)$ in finite time. Let $\phi$ be any maximal solution to \eqref{eq:NLsys}. Assume by contradiction that for all $j\in\dom\phi=\mathbb{N}$, $\phi(j)\notin L_V(c)$. Then, from \eqref{eq:DeltaVMain}
$$
V(\phi(j+1))-e^{\mu} V(\phi(j))\leq 0\quad\forall j\in\dom\phi
$$
which yields
\begin{equation}
\label{eq:DeltaVContr}
V(\phi(j))\leq e^{\mu j} V(\phi(0))\quad\forall j\in\dom\phi
\end{equation}
%Notice that if $V(\phi(0))\leq c$, from %\eqref{eq:DeltaVMain} one has that %$V(\phi(1))=0$, i.e., $\phi(1)\in S_u$. %Thus, since $S_u$ is forward invariant %for \eqref{eq:NLsys}, for all $j\geq 1$, %$\phi(j)\in S_u$. Hence, it is enough to %prove that for some $j^\star$, %$V(\phi(j^\star))\leq c$. Assume by %contradiction that for all %$j\in\dom\phi$, $V(\phi(j))>c$. Then, %from \eqref{eq:DeltaVMain},
%$$
%V(\phi(j+1))-V(\phi(j))\leq -c
%$$
%which yields
%\begin{equation}
%\label{eq:DeltaVContr}
%V(\phi(j))\leq -cj+V(\phi(0))\quad\forall %j\in\dom\phi
%\end{equation}
Pick 
$
\bar{j}=\left\lceil\frac{1}{\mu}\ln{\left(\frac{V(\phi(0))}{c}\right)}\right\rceil
$. 
Then, from \eqref{eq:DeltaVContr}, one gets 
$
V(\phi(\bar{j}))\leq c
$, 
which contradicts the fact that $V(\phi(j))>c$ for all $j\in\dom\phi$. In particular, define for all $x\in\R^n$
$$
\Gamma(x)\coloneqq \begin{cases}
\left\lceil\frac{1}{\mu}\ln{\left(\frac{V(x)}{c}\right)}\right\rceil&\text{if}\,\, x\in\R^n\setminus L_V(c)\\
0 &\text{else}  
\end{cases}
$$
which is locally bounded and nonnegative. The steps carried out so far show that for any maximal solution $\phi$ to \eqref{eq:NLsys}, 
$j\geq \Gamma(\phi(0))$ implies that $\phi(j)\in L_V(c)$. Thus, to conclude, let for all $r\geq 0$
$$
T(r)\coloneqq\sup_{x\in L_V(c)+r\mathbb{B}}\Gamma(x)
$$
notice that for all $r\geq 0$, $T(r)$ is finite, due to $\Gamma$ being locally bounded and $L_V(c)+r\mathbb{B}$ compact, and nonnegative due to $\Gamma$ being so. In particular, the definition of $T$ ensures that for any maximal solution $\phi$ to \eqref{eq:NLsys}, $j\geq T(\vert\phi(0)\vert_{L_V(c)})$ implies $\phi(j)\in L_V(c)$. This concludes the proof. \QEDA
\end{pf}
\subsection{Preliminaries on matrix-valued functions}
We consider matrix valued functions of the form:
\begin{equation}
\label{eq:functionX}
X\colon \mathcal{S}\to\mathcal{Y}
\end{equation}
where $\mathcal{S}$ is a finite dimensional real linear vector space and $\mathcal{Y}\subset\R^{n\times m}$. 
\medskip

\begin{definition}(Differential)
Let $X$ be defined as in \eqref{eq:functionX}. We say that $X$ is differentiable at $x\in\mathcal{S}$ if there exists a linear map $DX(x)\colon\mathcal{S}\to\mathcal{Y}$ such that:
$$
\lim_{\Vert h\Vert_{\mathcal{S}}\to 0}\frac{\Vert X(x+h)-X(x)-DX(x)h\Vert_{\mathcal{Y}}}{\Vert h\Vert_{\mathcal{S}}}=0
$$
where $\Vert\cdot\Vert_{\mathcal{S}}$ and $\Vert\cdot\Vert_{\mathcal{Y}}$ are any norms, respectively, on $\mathcal{S}$ and $\mathcal{Y}$.\hfill$\diamond$
\end{definition}
\medskip

\begin{definition}[\cite{shapiro1997first}]
Let $\mathcal{C}\subset\mathcal{S}$ be convex and $\mathbb{S}^n$ be the set of $n\times n$ symmetric matrices. A function $X\colon \mathcal{C}\to\mathbb{S}^n$ is said to be positive semidefinite convex (\emph{psd-convex}) on $\mathcal{C}$ if for all $x, y\in\mathcal{C}$ and $t\in [0, 1]$ the following holds:
$$
X(t x+(1-t)y)\preceq tX(x)+(1-t)X(y).
$$
Furthermore, we say that $X$ is positive semidefinite concave (\emph{psd-concave}) if $-X$ is psd-convex.\hfill$\diamond$
\end{definition}
\medskip

\begin{lemma}(\cite{dinh2011combining})
\label{lemm:DiffConv}
Let $\mathcal{C}\subset\mathcal{S}$ be convex, $\mathbb{S}^n$ be the set of $n\times n$ symmetric matrices, and $X\colon \mathcal{C}\to\mathbb{S}^n$ be differentiable on an open neighborhood of $\mathcal{C}$. Then, $X$ is psd-convex on $\mathcal{C}$ if and only if for all $x, y\in\mathcal{C}$
$$
X(y)-X(x)\succeq DX(x)(y-x). 
$$\hfill$\diamond$
\end{lemma}
\end{document}